\documentclass[twoside,leqno]{article}

\usepackage[letterpaper]{geometry}

\usepackage{siamproceedings}

\usepackage[T1]{fontenc}
\usepackage{amsfonts}
\usepackage{graphicx}
\usepackage{epstopdf}
\usepackage{enumitem}
\usepackage{algorithmic}
\ifpdf
  \DeclareGraphicsExtensions{.eps,.pdf,.png,.jpg}
\else
  \DeclareGraphicsExtensions{.eps}
\fi

\usepackage{color,latexsym,amsmath,amssymb, graphicx,revsymb, enumerate,xcolor,mathtools}
\renewcommand{\dim}{\mathrm{dim}\,}

\newcommand{\RR}{\mathbb{R}}
\renewcommand{\CC}{\mathbb{C}}

\newcommand{\FP}{\mathbb{F}_p}
\newcommand{\eps}{\varepsilon}
\newcommand{\cD}{\mathcal{D}}
\newcommand{\cE}{\mathcal{E}}

\newcommand{\tubes}{\mathbb{T}}

\newsiamremark{remark}{Remark}
\newsiamremark{hypothesis}{Hypothesis}
\crefname{hypothesis}{Hypothesis}{Hypotheses}
\newsiamthm{claim}{Claim}

\newsiamthm{conj}{Conjecture}

\newtheorem*{KMFC-B}{Conjecture \ref{kakeyaMaxmlConjA}$'$}
\newtheorem*{KMFC-C}{Conjecture \ref{kakeyaMaxmlConjA}$''$}

\usepackage{amsopn}

\begin{document}

\newcommand\relatedversion{}

\title{\Large A Survey of the Kakeya conjecture, 2000-2025}
    \author{Joshua Zahl\thanks{Chern Institute of Mathematics, Nankai University (\email{jzahl@nankai.edu.cn}).}}

\date{}

\maketitle

\begin{abstract} We survey progress on the Kakeya conjecture in Euclidean space.
\end{abstract}


\section{Introduction.}
A \emph{Besicovitch set} is a compact set $K\subset\RR^n$ that contains a unit line segment pointing in every direction. In \cite{Besicovitch2019}, Besicovitch constructed a Besicovitch set $K\subset\RR^2$ with Lebesgue measure 0. As a consequence, for every $\eta>0$ there exists $\delta>0$ and a set $\tubes$ of rectangles of dimensions $1\times\delta$ pointing in $\delta$-separated directions with the property that
\[
\sum_{\tubes}|T| \geq 1,\qquad \Big|\bigcup_{\tubes}T\Big| <\eta.
\]
We say that such a set $\tubes$ exhibits the \emph{Besicovitch compression phenomenon}. 

The \emph{Kakeya conjecture} is a family of closely related conjectures that seek to understand and quantify the Besicovitch compression phenomenon.

\begin{conj}[Kakeya set conjecture]\label{kakeyaSetConj}
Every Besicovitch set $K\subset\RR^n$ has Minkowski and Hausdorff dimension $n$.
\end{conj}

Conjecture \ref{kakeyaSetConj} is typically studied by first discretizing the set $K$ at a small scale $\delta>0$. $N_\delta(K)$ contains a $\delta$ tube (the $\delta$-neighbourhood of a unit line segment) pointing in every direction, and in particular $N_\delta(K)$ contains a $\delta$ tube pointing in each $\delta$-separated direction. This gives rise to a discretized analogue of Conjecture \ref{kakeyaSetConj}, which we state as Conjecture \ref{KakeyaDiscretized} below. Before stating the conjecture, we introduce the following notation: if $f$ and $g$ are functions of $\delta$, we say $f\lessapprox g$ if for all $\eps>0$, there exists $C_\eps$ (which may also depend on the ambient dimension $n$) so that $f\leq C_\eps\delta^{-\eps} g$. For example $\log 1/\delta\lessapprox 1$. For the purpose of this survey we will sometimes be intentionally vague when distinguishing between inequalities of the form $A \lesssim B$ (or $A = O(B)$) and $A\lessapprox B$. In contrast to delicate questions such as the $L^p$ boundedness of operators related to the Hilbert transform, the distinction between $A = O(B)$ and $A\lessapprox B$ is usually unimportant for Kakeya-type problems. 
\begin{conj}\label{KakeyaDiscretized}
(A): Let $\tubes$ be a set of $\delta$ tubes in $\RR^n$ pointing in $\delta$-separated directions. Then
\[
\Big|\bigcup_{\tubes}T\Big| \gtrapprox \sum_\tubes |T|.
\]

\medskip
(B): For each $T\in\tubes$, let $Y(T)\subset T$ be measurable with $|Y(T)|\geq (\log 1/\delta)^{-1} |T|$. Then
\[
\Big|\bigcup_{\tubes}T\Big| \gtrapprox \sum_\tubes |T|.
\]
\end{conj}
Conjecture \ref{KakeyaDiscretized}(A) implies that every Besicovitch set in $\RR^n$ has Minkowski dimension $n$, while Conjecture \ref{KakeyaDiscretized}(B) implies that every Besicovitch set in $\RR^n$ has Hausdorff dimension $n$. 

The following variant of the Kakeya conjecture quantifies the size of the region in which the tubes from $\tubes$ may substantially overlap:

\begin{conj}[Kakeya maximal function conjecture]\label{kakeyaMaxmlConjA}
Let $\tubes$ be a set of $\delta$ tubes in $\RR^n$ pointing in $\delta$-separated directions. Then
\begin{equation}\label{kakeyaMaximalFnConj}
\Big\Vert \sum_{\tubes}\chi_T \Big\Vert_{p}\lessapprox 1,\quad 1\leq p \leq \frac{n}{n-1}.
\end{equation}
\end{conj}

Conjecture \ref{kakeyaMaxmlConjA} is called the \emph{Kakeya maximal function conjecture}. When $p=1$, \eqref{kakeyaMaximalFnConj} follows from the fact that $\#\tubes\leq \delta^{1-n}$. By considering the example where $\tubes$ consists of $\delta^{1-n}$ tubes that contain the origin, it is straightforward to verify that \eqref{kakeyaMaximalFnConj} is false for $p>n/(n-1)$. It is sometimes more convenient to write Conjecture \ref{kakeyaMaxmlConjA} in the following equivalent form:
\begin{KMFC-B}
Let $\tubes$ be a set of $\delta$ tubes in $\RR^n$ pointing in $\delta$-separated directions. Let $\lambda\in[\delta,1]$ and for each $T\in\tubes,$ let $Y(T)\subset T$ be measurable with $|Y(T)|\geq \lambda|T|$. Then
\[
\Big|\bigcup_{\tubes}Y(T)\Big| \lessapprox \lambda^n \sum_{T\in\tubes}|T|.
\]
\end{KMFC-B}

Finally, as the name suggests, we can express Conjecture \ref{kakeyaMaxmlConjA} as a statement about the $L^p$ boundedness of the \emph{Kakeya maximal operator}: for $\delta>0$, $f\colon\RR^n\to\CC$ and $e\in S^{n-1}$, define the Kakeya maximal operator 
\[
\mathcal{K}_\delta f(e) = \sup_{T \| e} \frac{1}{|T|}\int_T |f|,
\]
where the supremum is taken over all $\delta$ tubes $T$ pointing in the direction $e$.

\begin{KMFC-C}
\[
\Vert \mathcal{K}_\delta f\Vert_{L^n(S^{n-1})}\lessapprox  \Vert f\Vert_{L^n(\RR^n)}.
\]
\end{KMFC-C}

Conjectures \ref{kakeyaMaxmlConjA}, \ref{kakeyaMaxmlConjA}$'$ and \ref{kakeyaMaxmlConjA}$''$ are equivalent. They imply Conjecture \ref{KakeyaDiscretized}(B), which in turn implies Conjecture \ref{KakeyaDiscretized}(A).

The Kakeya conjecture is closely connected to several open problems in Fourier analysis. A seminal paper by Bourgain \cite{Bourgain1991} from 1991 explored the connection between the Kakeya maximal function conjecture and Stein's Fourier restriction conjecture. For $f\colon \RR^{n-1}\to\CC$ and $x\in\RR^n$, we define the \emph{Fourier extension operator} 
\[
Ef(x) = \int_{B(0,1)} f(\omega)e^{2\pi i(x_1\omega_1 + \ldots + x_{n-1}\omega_{n-1} + x_n|\omega|^2)}d\omega.
\]

In its adjoint form, Stein's Fourier restriction conjecture asserts that 
\[
\Vert Ef\Vert_{L^p(\RR^n)}\leq C_{n,p}\Vert f\Vert_{L^p(\RR^{n-1})},\qquad p>\frac{2n}{n-1}.
\]
The conjecture was proved by Fefferman and Stein \cite{Fefferman1970} for $n=2$, and remains open in dimension $n\geq 3$. Bourgain analyzed the extension operator $Ef$ by first localizing to a large ball $B(0,R)$, and then decomposing the function $f$ into a sum $f = \sum f_i$, where each $f_i$ is supported on a disk $\theta_i$ of radius $R^{-1/2}$. The functions $Ef_i$ can be decomposed into a sum of wave packets $Ef_i = \sum_j g_{i,j}$, where each $g_{i,j}$ is essentially supported on a long (relatively) thin tube whose direction is determined by the location of the disk $\theta_i$. Thus in order to analyze the extension operator $Ef$, it suffices to understand the possible intersection patterns of families of tubes in $\RR^n$ pointing in different directions (this is closely related to the Kakeya conjecture), plus the possible patterns of constructive and destructive interference between the functions supported on these tubes.

The goal of this article is to survey progress on Conjectures \ref{kakeyaSetConj} and \ref{kakeyaMaxmlConjA}, with an emphasis on developments that occurred after the excellent 2002 survey paper \cite{KatzTao2002}. In particular, we will not discuss the techniques and results of Davies \cite{Davies1971} and Cordoba \cite{Cordoba1977}, who proved Conjectures \ref{kakeyaSetConj} and \ref{kakeyaMaxmlConjA} in $\RR^2$. Nor will we discuss the geometric arguments of Drury \cite{Drury1983}; Christ-Duoandikoetxea-Rubio de Francia \cite{ChristDuoandikoetxeaFrancia1986}; and Bourgain \cite{Bourgain1991}, or the arguments from additive combinatorics introduced by Bourgain \cite{Bourgain1999} and Katz-Tao \cite{KatzTao1999, KatzTao2002b}. We refer the reader to the previous surveys \cite{Wolff1999, KatzTao2002} for a discussion of these earlier developments, and to Chapter 10 of Wolff's excellent book \cite{ Wolff2003} for an introduction and further background on the Kakeya problem. Finally, we will restrict attention to the Kakeya conjecture in Euclidean space; in particular we will not discuss Dvir's solution to the finite field Kakeya conjecture \cite{Dvir2009}, or the subsequent developments \cite{SarafSudan2008, DvirKoppartySarafSudan2012, BukhChao2021, Arsovski2024, Dhar2024} in the finite field and $p$-adic setting.


\section{Kakeya in $\RR^3$.}
In \cite{Wolff1995}, Wolff proved a Kakeya maximal function estimate in $\RR^n$ at dimension $\frac{n+2}{2}$.  Wolff's argument not only applies to sets of $\delta$ tubes pointing in $\delta$-separated directions, but also applies to sets of tubes that satisfy the following non-clustering condition:

\begin{definition}\label{defnWolffAxioms}
Let $\tubes$ be a set of $\delta$ tubes in $\RR^n$. We say that $\tubes$ satisfies the \emph{Convex Wolff Axioms} if every convex set $W\subset\RR^n$ contains $O(|W|\delta^{1-n})$ tubes from $\tubes$. 
\end{definition}
In the discussion below it will sometimes be helpful to weaken Definition \ref{defnWolffAxioms} and instead require that each convex set $W$ contains $\lessapprox |W|\delta^{1-n}$ tubes from $\tubes$. We also remark that \cite{WangZahl2025} uses the terminology \emph{Katz-Tao} Convex Wolff Axioms for this definition; we will use the terms ``Convex Wolf Axioms'' and ``Katz-Tao Convex Wolff Axioms'' interchangeably.

\begin{theorem}[Wolff \cite{Wolff1995}]\label{WolffThm}
Let $n\geq 2$ and let $\tubes$ be a set of $\delta$ tubes in $\RR^n$ that satisfy the Convex Wolff Axioms. Then
\begin{equation}
\Big|\bigcup_{\tubes}T\Big| \gtrapprox \delta^{\frac{n-2}{2}}\sum_{\tubes}|T|,
\end{equation}
and
\begin{equation}
\Big\Vert\sum_{\tubes}\chi_T\Big\Vert_{\frac{d}{d-1}} \lessapprox \Big(\frac{1}{\delta}\Big)^{\frac{n}{d}-1},\quad d = \frac{n+2}{2}.
\end{equation}
\end{theorem}

In $\RR^3$, Theorem \ref{WolffThm} says that every Besicovitch set has Hausdorff dimension at least $5/2$, and every union of $\delta$ tubes of cardinality $\delta^{-2}$ satisfying the Convex Wolff Axioms has volume $\big|\bigcup_\tubes T\big|\gtrapprox\delta^{1/2}$. Wolff's argument only uses simple facts about the volume of intersections of tubes and other thickened neighbourhoods of affine subspaces of $\RR^n$. As we will discuss below, in $\RR^3$ Theorem \ref{WolffThm} is best-possible if one only uses these tools.


\subsection{Near misses to the Kakeya conjecture in $\RR^3$: the Heisenberg group and $SL_2$ examples.}\label{HeisenbergSL2ExampleSection}
A \emph{near miss} to the Kakeya set conjecture is a set that has properties in common with a Besicovitch set, but which has dimension (or some suitable analogue, such as cardinality) smaller than what is predicted by the Kakeya set conjecture. A near miss to the Kakeya set conjecture helps us understand the limitations of certain proof techniques---if a proof technique cannot distinguish between a genuine Kakeya set in $\RR^n$ and a near-miss, then this proof technique cannot prove the Kakeya set conjecture. A near miss to the Kakeya set conjecture thus forms a barrier to proving the Kakeya set conjecture. This is (loosely) analogous to the relativization, natural proofs, and algebrization barriers from complexity theory. 

In \cite{KatzLabaTao2000}, Katz, \L{}aba, and Tao discovered the Heisenberg group example:
\begin{equation}
\mathbb{H} = \{(x,y,z)\in\CC^3\colon \operatorname{Im}(z) = \operatorname{Im}(x\bar y)\}.
\end{equation}
For $a,b,c,d\in\RR$, we can verify that the line 
\[
L = \big\{(0, c+di,a)+s(1, b, c-di)s\colon s \in \CC\big\}
\] 
is contained in $\mathbb{H}$. Let $\rho = \delta^{1/2}$ and let $E$ be the closure of $N_\rho(\mathbb{H})\cap B(0,1)$. $E$ is a compact subset of $\CC^3$ that contains a set of $\delta^{-2}$ (complex) tubes (i.e.~sets $T\subset\CC^3$ that are the thickened neighbourhood of a (complex) unit line segment, with $|T| \sim\delta^{2}$) that satisfy the Convex Wolff Axioms. On the other hand, we have $|E|\sim\delta^{1/2}$. We say that the Heisenberg group example is a near miss to the Kakeya set conjecture in three dimensions, at dimension $5/2$. 

In \cite{KatzLabaTao2000}, Katz, \L{}aba, and Tao proved that every Besicovitch set in $\RR^3$ has upper Minkowski dimension at least $5/2+c_0$ for a (small) absolute constant $c_0$; we will discuss their proof below. They did this by exploiting the fact that a Besicovitch set contains a line pointing in every direction in $\RR^3$. The (complex) lines in the Heisenberg group example, on the other hand, point in a positive codimension set of directions in $\CC^3$.

In \cite{KatzZahl2019}, Katz and the author discovered a second near-miss to the Kakeya conjecture at dimension $5/2$. They named this the $SL_2$ example. Let $R$ be the ring $\FP[t]/(t^2)$.  Define
\begin{equation}\label{defnOfSL2}
X = \{(x_1+x_2t,\ y_1+y_2t,\ z_1+z_2t)\in R^3\colon z_2 = x_1y_2 - x_2y_1\}.
\end{equation}
We have that $\# X= p^5 = (\#R)^{5/2}$. We can verify that for $a,b,c,d,e\in \FP$ with $ad-bc=1$, the line
\[
L = \big\{(a+eat, b+ebt, 0) + (c+ect, d+edt, 1)s \colon s\in R \big\}
\]
is contained in $X$. Thus $X$ contains a set of roughly $p^4 = (\#R)^2$ lines, and these lines satisfy an analogue of the Convex Wolff Axioms (specifically, every plane contains $O(\#R)$ lines).  

The intuition for this example is as follows. Each number $x\in[0,1]\subset\RR$ can be written as $x = \delta^{1/2}x_1 + \delta x_2+O(\delta)$, where $x_1$ and $x_2$ are integers between $0$ and $\lfloor\delta^{-1/2}\rfloor$. The ring $R$ is meant to model this two-scale decomposition of the interval $[0,1]$. Elements of $\FP\subset \FP[t]/(t^2)$ represent the coarse (i.e.~$\delta^{1/2}$) scale, while elements of $t\cdot\FP$ represent the fine (i.e.~$\delta$) scale. In particular, the $SL_2$ example resembles a Besicovitch set in $\RR^3$ that has Hausdorff dimension $5/2$, but upper Minkowski dimension $3$. It is for this reason that Katz, \L{}aba, and Tao did not encounter the $SL_2$ example near-miss when proving that every Besicovitch set in $\RR^3$ has upper Minkowski dimension $5/2+c_0$. 


\subsection{Upper Minkowski dimension $5/2+c_0$: Overcoming the Heisenberg group example.}\label{52KLTSection}

\begin{theorem}[Katz-\L{}aba-Tao \cite{KatzLabaTao2000}]\label{KLTTheorem}
Let $\tubes$ be a set of $\delta$ tubes in $\RR^3$ of cardinality $\delta^{-2}$ pointing in $\delta$-separated directions. Then for at least one scale $\rho\in\{\delta, \delta^{1/2},\delta^{1/4}\}$ we have
\begin{equation}\label{biggerAtSomeScale}
\Big|N_\rho\Big(\bigcup_{\tubes}T\Big)\Big|\gtrsim\rho^{1/2-c_0},
\end{equation}
where $c_0>0$ is an absolute constant. In particular, every Besicovitch set in $\RR^3$ has upper Minkowski dimension at least $5/2+c_0$.
\end{theorem}

As discussed above, Katz, \L{}aba, and Tao's proof of Theorem \ref{KLTTheorem} used the fact that a Besicovitch set in $\RR^3$ contains a line pointing in every direction---this is a property that does not hold for the Heisenberg group example. To make use of this property, Katz, \L{}aba, and Tao used the following ``sums versus differences'' result from additive combinatorics:
\begin{theorem}\label{sumDifferences}
Let $G$ be an Abelian group, and let $A,B\subset G$ with $\#A\leq N,$ $\#B\leq N,$ and $\#(A+B)\leq N$. Then $\#(A-B)\lesssim N^{2-c}$, where $c>0$ is an absolute constant.
\end{theorem}
A weaker form of Theorem \ref{sumDifferences} with $o(N^2)$ in place of $N^{2-c}$ was first proved by Balog and Szemer\'edi \cite{BalogSzemeredi1994}, and the version stated above follows from the work of Gowers \cite{Gowers1998}. Theorem \ref{sumDifferences} as stated was proved and first used in the context of the Kakeya problem by Bourgain \cite{Bourgain1999}. 

The idea is to apply Theorem \ref{sumDifferences} to three slices of the Besicovitch set---the sets $A$, $B$, and $(A+B)/2$ correspond to the intersection of $\bigcup_\tubes T$ with the planes $\{x=0\}$, $\{x=1\}$, and $\{x=1/2\}$ respectively, while $A-B$ represents the set of directions of the tubes in $\tubes$. However, a naive application of this method only establishes that Besicovitch sets in $\RR^3$ have dimension at least $2+c$ for a small constant $c>0$. 

To prove Theorem \ref{KLTTheorem}, Katz, \L{}aba and Tao analyzed the structure of a (hypothetical) Besicovitch set in $\RR^3$ with Minkowski dimension exactly $5/2$. If such a Besicovitch set is discretized at scale $\delta$, it yields a set $\tubes$ of $\delta$ tubes pointing in $\delta$-separated directions with the following property: $\#\tubes \sim \delta^{-2}$, and for both $\rho = \delta^{1/2}$ and $\rho = \delta^{1/4}$, there exists a set $\tubes_\rho$ of $\rho$ tubes pointing in $\rho$-separated directions with the property that each $T\in\tubes$ is contained in a $\rho$ tube, and each $\rho$ tube contains about $(\rho/\delta)^2$ tubes from $\tubes$. This motivates the following definition. In what follows, we say that a set of $\rho$ tubes $T_\rho$ (or more generally, a collection of convex sets) is a \emph{cover}\footnote{more precisely, an almost-balanced, almost-partitioning cover} of $\tubes$ if each tube in $\tubes$ is contained in at least one, and at most $O(1)$ tubes from $\tubes_\rho$, and each tube from $\tubes_\rho$ contains approximately the same number of tubes from $\tubes$.

\begin{definition}\label{defnSticky}
Let $\tubes$ be a set of $\delta$ tubes in $\RR^n$ of cardinality $\#\tubes\sim \delta^{1-n}$. We say that $\tubes$ is \emph{sticky} if for all $\delta\leq\rho\leq 1$, there exists a cover $\tubes_\rho$ of $\tubes$ that satisfies the Convex Wolff Axioms.
\end{definition}
Stickiness was first introduced in \cite{KatzLabaTao2000}, though our definition above is motivated by more recent results and is slightly anachronistic.

Next, Katz, \L{}aba, and Tao proved that the set of tubes $\tubes$ described above must have two structural properties that they named \emph{planiness} and \emph{graininess}. Planiness says that for a typical point $x\in\bigcup_\tubes T$, most of the tubes containing $x$ make angle $\leq\theta$ with a common affine plane $H_x$, for some $\theta<\!\!<1$. We now know, as a consequence of the multilinear Kakeya theorem (see Theorem \ref{multilinearKakeyaThm} below), that this property holds for every set of $\delta$ tubes in $\RR^3$ for which $\big|\bigcup_\tubes T\big|$ is substantially smaller than $\sum_\tubes|T|$. Graininess says that $\bigcup_\tubes T$ can be decomposed as a union of mostly-disjoint rectangular prisms, called grains, of dimensions $\delta\times\delta^{1/2}\times\delta^{1/2}$; nearby prisms have almost parallel normal vectors, in the sense that the normal vectors of rectangular prisms inside a common $\delta^{1/2}$ ball agree up to resolution $\delta^{1/2}$. Guth \cite{Guth2016b} later proved that a similar property holds for every set of $\delta$ tubes in $\RR^3$ (satisfying mild hypotheses) for which $\big|\bigcup_\tubes T\big|$ is substantially smaller than $\sum_\tubes|T|$, though the size of the corresponding grains becomes smaller as $\big|\bigcup_\tubes T\big|$ approaches $\sum_\tubes|T|$.

Finally, after showing that a (hypothetical) Besicovitch set in $\RR^3$ with upper Minkowski dimension $5/2$ gives rise to a set of $\delta$ tubes that are sticky, plany, and grainy, Katz, \L{}aba, and Tao completed their proof of Theorem \ref{KLTTheorem} by constructing an Abelian group $G$ and sets $A,B\subset G$ that violate the conclusion of Theorem \ref{sumDifferences}. 

We conclude this section by briefly discussing extensions of the above ideas to the Kakeya problem in higher dimensions. \L{}aba and Tao \cite{LabaTao2001} showed that every Besicovitch set in $\RR^n$ with upper Minkowski dimension close to $\frac{n+2}{2}$ must be sticky, plany, and grainy. They used this to establish an analogue of Theorem \ref{KLTTheorem} for $n\geq 4$. Bourgain \cite{Bourgain1999} and subsequently Katz-Tao \cite{KatzTao1999, KatzTao2002b} applied Theorem \ref{KLTTheorem} and the ``three slices'' idea described above (later extended to a larger number of slices) to obtain new estimates in the direction of Conjectures \ref{kakeyaSetConj} and \ref{kakeyaMaxmlConjA}.


\subsection{Hausdorff dimension $5/2+c_1$: Overcoming the $SL_2$ example.}
\begin{theorem}[Katz-Zahl \cite{KatzZahl2019}]\label{KZTheorem}
Let $\tubes$ be a set of $\delta$ tubes in $\RR^3$ of cardinality $\delta^{-2}$ that satisfy the Convex Wolff Axioms. For each $T\in\tubes,$ let $Y(T)\subset T$ with $|Y(T)|\geq(\log 1/\delta)^{-1}|T|$. Then
\begin{equation}\label{conclusionKatzZahl}
\Big|\bigcup_{\tubes}Y(T)\Big|\gtrsim \delta^{1/2-c_1},
\end{equation}
where $c_1>0$ is an absolute constant. In particular, every Besicovitch set in $\RR^3$ has Hausdorff dimension at least $5/2+c_1$. 
\end{theorem}
Rather than requiring that the tubes in $\tubes$ point in $\delta$-separated directions, Theorem \ref{KZTheorem} only requires that the tubes satisfy the Convex Wolff Axioms. Katz and the author's proof of Theorem \ref{KZTheorem} makes crucial use of the fact that the tubes in $\RR$ are subsets of real Euclidean space---this is a property that does not hold for the Heisenberg group example. To make use of this property, Katz, and the author used the following discretized projection theorem due to Bourgain \cite{Bourgain2010}:

\begin{theorem}\label{bourgainThm}
Let $E\subset [0,1]^2$ be a union of $\delta$ balls, with $|E| \sim \delta$. Suppose that $E$ satisfies the Frostman non-concentration condition
\[
|E\cap B(x,r)|\lessapprox r|E|\quad\textrm{for all balls}\ B(x,r),\ \delta\leq r\leq 1.
\]
Let $\Lambda\subset [0,1]$ be a union of $\delta$-intervals satisfying the Frostman non-concentration condition
\[
|\Lambda\cap I|\lessapprox |I|^{1/100}|\Lambda|\quad\textrm{for all intervals}\ I\subset S^1.
\]
Then there exists $\lambda\in\Lambda$ so that
\[
|E\cdot(1,\lambda)|\gtrapprox |E|^{1/2-c_2},
\]
where $c_2>0$ is an absolute constant. 
\end{theorem}

The analogue of Theorem \ref{bourgainThm} is false over $\CC$: define $E_{\RR}$ to be the $\delta$-neighbourhood of $(\RR\times\RR)\cap B(0,1)\subset\CC^2$ and define $\Lambda_{\RR} = [0,1]\subset\RR\subset\CC$. Then $E_{\RR}$ and $\Lambda_{\RR}$ satisfy the hypotheses of Theorem \ref{bourgainThm} (suitably modified), but we have $|E_{\RR}\cdot(1,\lambda)|\sim |E_{\RR}|^{1/2}$ for all $\lambda\in\Lambda_{\RR}$. When proving Theorem \ref{KZTheorem}, there is a key step where Theorem \ref{bourgainThm} is used. If one instead applies the proof of Theorem \ref{KZTheorem} step-by-step to the Heisenberg group example $\mathbb{H}$, then at the step where Theorem \ref{bourgainThm} is used, the inputs $E$ and $\Lambda$ to the theorem are (closely related to) $E_{\RR}$ and $\Lambda_{\RR}$. In particular, the conclusion of Theorem \ref{bourgainThm} is false at this key step. 

Katz and the author's proof of Theorem \ref{KZTheorem} loosely parallels an analogous proof in the finite field setting due to Bourgain, Katz, and Tao \cite{BourgainKatzTao2004} (this latter proof uses a finite field analogue of the Erd\H{o}s-Szemer\'edi sum-product theorem \cite{ErdosSzemeredi1983} in place of Theorem \ref{bourgainThm}). In brief, the proof is as follows. Suppose for contradiction that there is a set $\tubes$ of $\delta$ satisfying the hypotheses of Theorem \ref{bourgainThm} for which \eqref{conclusionKatzZahl} fails. After pigeonholing it is possible to find pairs of tubes $T_1,T_2$ and $T_A,T_B$ and a set $\tubes_{1,2}\subset\tubes$ (resp.~$\tubes_{A,B}\subset\tubes$) of cardinality roughly $\delta^{-1}$ that intersect both $T_1$ and $T_2$ (resp.~$T_A$ and $T_B$). Each tube in $T_{1,2}$ intersects roughly $\delta^{-1/2}$ tubes from $\tubes_{A,B}$, and vice-versa. The situation thus far is consistent with both the Heisenberg group and $SL_2$ examples. If one applies these arguments to the Heisenberg group example, then for a typical tube $T\in\tubes_{1,2}$, the points of intersection from tubes in $T'\in\tubes_{A,B}$, i.e. the set $\{T'\cap T\colon T'\in\tubes_{A.B}\}$ are spread-out along $T$. In fact they obey a Frostman-type non-concentration condition. At the opposite extreme, if one applies these arguments to the $SL_2$ example, then the set $\{T'\cap T\colon T'\in\tubes_{A.B}\}$ is concentrated into a short tube segment (the analogue of a tube segment of length $\delta^{1/2}$). 

Guided by the behavior of these two examples, Katz and the author showed that every counter-example to Theorem \ref{KZTheorem} must be either of ``Heisenberg type'' or ``$SL_2$ type.''  Different arguments are used to handle these two cases.

If $\tubes$ is of Heisenberg type, then we can identify each tube $T\in\tubes_{A,B}$ with a $\delta$ ball in $\RR^2$, and we can identify the tubes in $\tubes_{1,2}$ with $\delta$ intervals in $[0,1]$ (this latter identification is not one-to-one), so that the corresponding union of points and intervals violates Theorem \ref{bourgainThm}. This is impossible, and we conclude that the original set $\tubes$ must have satisfied \eqref{conclusionKatzZahl}. 

If $\tubes$ is of $SL_2$ type, then elementary geometric arguments are used to show that the tubes in $\tubes$ must cluster into twisting planks called \emph{regulus strips}, and these regulus strips must intersect in a carefully prescribed manner that violates (a curved variant of) the Szemer\'edi-Trotter theorem. This is impossible, and we (again) conclude that the original set $\tubes$ must have satisfied \eqref{conclusionKatzZahl}. 

While the ideas used to prove Theorem \ref{KZTheorem} did not play a substantial role in the resolution of Conjecture \ref{kakeyaSetConj} in $\RR^3$, the analysis of unions of regulus strips continues to be an area of active study that is closely related to questions in projection theory in Euclidean space and the first Heisenberg group.


\section{Hausdorff dimension $3$: the Wang-Zahl proof.}
In the trilogy of papers \cite{WangZahl2022, WangZahl2024, WangZahl2025}, Wang and the author proved the Kakeya set conjecture in $\RR^3$. The precise statement is as follows.
\begin{theorem}\label{WangZahlWolffAxiomThm}
Let $\tubes$ be a set of $\delta$ tubes in $\RR^3$ that satisfy the Convex Wolff Axioms. Then
\[
\Big|\bigcup_{T\in\tubes}T\Big| \gtrapprox \sum_{T\in\tubes}|T|.
\]
As a consequence, Conjectures \ref{kakeyaSetConj} and \ref{KakeyaDiscretized} are true for $n=3$.
\end{theorem}
In particular, other than the Heisenberg group and $SL_2$ examples, there do not appear to be any additional near misses to the Kakeya set conjecture in $\RR^3$ (though there may still be as-yet undiscovered near misses to the Kakeya maximal function conjecture in $\RR^3$). A recent article by Guth \cite{Guth2025} does an excellent job both explaining the key ideas of the proof of Theorem \ref{WangZahlWolffAxiomThm}, and also presenting substantial new and innovative simplifications to the core of the argument. As such, we will be brief in our remarks here and encourage the reader to consult \cite{Guth2025} (and after that \cite{WangZahl2022, WangZahl2024, WangZahl2025}), for further details. 


\subsection{The sticky case.}

In the early 2000s, Katz and Tao formulated a program for proving the Kakeya set conjecture in $\RR^3$. They hypothesized that every counter-example to Conjecture \ref{KakeyaDiscretized}(A) in $\RR^3$ must have the structural properties stickiness, planiness, and graininess described in Section \ref{52KLTSection}. They then formulated a plan to leverage these properties to show that a counter-example to Conjecture \ref{KakeyaDiscretized}(A) would imply the existence of a set that violates Bourgain's discretized sum-product theorem \cite{Bourgain2003} which was a precursor to (and close relative of) Theorem \ref{bourgainThm}. In a talk and accompanying blog post \cite{Tao2014}, Tao explained the Katz-Tao program for solving the Kakeya set conjecture in $\mathbb{R}^3$. At the time, it was not clear whether every (hypothetical) counter-example to Conjecture \ref{KakeyaDiscretized}(A) must be sticky. In particular, the discovery of the $SL_2$ example suggested that in order to prove the Kakeya set conjecture in $\RR^3$, it would be necessary to study sets that are far from being sticky. 

In \cite{WangZahl2022}, Wang and the author side-stepped this issue by restricting attention to sticky sets of tubes $\tubes$ in $\RR^3$, i.e.~sets satisfying (a variant of) Definition \ref{defnSticky}. 

\begin{theorem}[Wang-Zahl \cite{WangZahl2022}]\label{stickyKakeyaThm}
Let $\tubes$ be a sticky set of $\delta$ tubes in $\RR^3$ that point in $\delta$-separated directions, with $\#\tubes\sim\delta^{-2}$. For each $T\in\tubes$, let $Y(T)\subset T$ with $|Y(T)|\geq(\log 1/\delta)^{-1}|T|$. Then
\[
\Big|\bigcup_{\tubes}Y(T)\Big|\gtrapprox 1.
\]
\end{theorem}
We remark that while the statement of Theorem \ref{stickyKakeyaThm} requires that the tubes point in $\delta$-separated directions, the proof does not use this fact in an essential way. In follow-up work (see Proposition \ref{FrostmanWolffEveryScale}) this hypothesis is removed.

We will briefly describe the proof. For simplicity, we will suppose that $Y(T) = T$ for each tube $T\in\tubes$. The proof follows the initial steps of the argument outlined in \cite{Tao2014} to show that a (hypothetical) counter-example to Theorem \ref{stickyKakeyaThm} must have certain structural properties analogous to those in the Heisenberg group example. We begin with the following observation of Katz and Tao: suppose that Theorem \ref{stickyKakeyaThm} is false, and let $\tubes$ be the ``worst'' counter-example, in the sense that $\big|\bigcup_\tubes T\big|=\delta^\sigma$ for some $\sigma>0$, and for all $\rho>0$ and all sticky collections $\tubes_\rho$ of $\rho$-tubes, we have $\big|\bigcup_{\tubes\rho}T_\rho\big|\geq\rho^{\sigma}$. Since $\sum_{\tubes}|T|\sim 1$, we have that on average, about $\delta^{-\sigma}$ tubes intersect a typical point of $\bigcup_\tubes T$, i.e.~the tubes in $\tubes$ have typical intersection multiplicity $\delta^{-\sigma}$; we will denote this quantity by $\mu_\tubes$.

Since $\tubes$ is sticky, for every intermediate scale $\rho\in[\delta,1]$ there exists a cover $\tubes_\rho$ of $\tubes$ that satisfies the Convex Wolff Axioms. After a refinement, we can show that $\tubes_\rho$ is also sticky, and thus $\big|\bigcup_{\tubes_\rho}T_\rho\big|\geq\rho^\sigma$, or equivalently, $\mu_{\tubes_\rho}\lesssim \rho^{-\sigma}$. Furthermore for each $T_\rho\in\tubes_\rho$, there is an (anisotropic) rescaling so that the tubes in the set $\tubes[T_\rho]:=\{T\in\tubes\colon T\subset T_\rho\}$ become $\delta/\rho$ tubes inside the unit ball. Define $\tubes^{T_\rho}$ to be the image of the tubes in $\tubes[T_\rho]$ under this rescaling. Then $\tubes^{T_\rho}$ forms a sticky collection of $\delta/\rho$ tubes, and hence $\mu_{\tubes[T_\rho]}\lesssim (\delta/\rho)^{-\sigma}$. We now have the key observation
\begin{equation}\label{muControl}
\delta^{-\sigma} \sim \mu_\tubes \leq \big(\mu_{\tubes_\rho}\big)\big(\sup_{T\in \tubes_\rho} \mu_{\tubes[T_\rho]}\big)\lesssim \rho^{-\sigma}(\delta/\rho)^{-\sigma}=\delta^{-\sigma}.
\end{equation}

The first inequality follows from the observation that the number of $\delta$ tubes intersecting a point $x\in\bigcup_{\tubes}T$ is at most the number of $\rho$ tubes from $\tubes_\rho$ that intersect $x$, times the maximum intersection multiplicity of $\delta$ tubes inside a $\rho$ tube. But since the LHS and RHS of Inequality \eqref{muControl} are equal, we have $\mu_{\tubes_\rho}\sim\rho^{-\sigma}$ and $\mu_{\tubes[T_\rho]}\sim (\delta/\rho)^{-\sigma}$ for each $\rho$ tube in $\tubes_\rho$ (when performing this proof rigorously, a similar but slightly weaker statement is true). This means that both $\tubes_\rho$ and the sets $\tubes^{T_\rho}$ are also ``worst'' counter-examples to Theorem \ref{stickyKakeyaThm}. We conclude that $\tubes$ is coarsely (or statistically) self-similar: every structural property of $\tubes$ at scale $\delta$ is mirrored at all locations and scales. 

A consequence of the multilinear Kakeya theorem plus coarse self-similarity is that $\bigcup_{\tubes}T$ is grainy. The locations of the grains at scale $\delta$ can be encoded as a (discretized) one-dimensional set $F\subset\RR^2$. The locations of the grains at the coarser scale $\delta^{1/2}$ can be encoded as a (discretized) one-dimensional set $G\subset\RR^2$. Examining the interplay between grains at these two scales, it is possible to find many closed paths inside $\bigcup_\tubes T$ that encode arithmetic identities between the sets $F$ and $G$: we conclude that there is an discretized (almost) Ahlfors-regular set $E\subset[0,1]$ of dimension $\dim(E)<1$ so that most of $(F-F)\cdot G$ is contained in $E$. 

The set $(F-F)\cdot G$ contains a scaled copy of the orthogonal projection of $G\subset\RR^2$ onto each one-dimensional subspace spanned by vectors of the form $e_{x_1,x_2}:=(x_1-x_2)/|x_1-x_2|$ for $x_1,x_2\in F$. In particular, for most pairs $x_1,x_2\in F$, $e_{x_1,x_2}\cdot G$ is contained in a discretized (almost) Ahlfors regular set of dimension strictly smaller than 1. Recent tools in projection theory \cite{OrponenShmerkinWang2024} (which at their core use Theorem \ref{bourgainThm}) show that this is impossible, unless $F$ and $G$ are contained in orthogonal lines. This is the key step where we use the fact that the tubes in $\tubes$ are subsets of $\RR^3$ rather than $\CC^3$, and thus we distinguish $\tubes$ from the Heisenberg group example: Over $\CC$, the sets $F = G =\RR^2\cap B(0,1)\subset\CC^2$ are not contained in orthogonal lines, but $e_{x_1,x_2}\cdot G\subset \RR$ (this latter set is Ahlfors regular and has dimension $\dim(\RR) = 1 < 2 = \dim(\CC)$) for every $x_1,x_2\in F$.

Finally, if $F$ and $G$ are contained in orthogonal lines, then this forces $\bigcup_\tubes T$ to (almost) look like a Cartesian product $[0,1]\times W$. More precisely, $\bigcup_\tubes T$ is the image of a Cartesian product $[0,1]\times W$ under a map that is linear on each $\{z=z_0\}$ slice, and twists (with controlled first and second derivatives) as $z$ varies. This allows us to project $\bigcup_{\tubes}T$ to a plane in such a way that the fibers of the projection are unit intervals. The image of each tube in $\tubes$ becomes the $\delta$ neighbourhood of a $C^2$ curve. To finish the proof we apply a variant \cite{PramanikYangZahl2022} of Wolff's circular maximal function estimate \cite{Wolff1997} to conclude that the union of these thickened plane curves has volume $\approx 1$, and hence (since each fiber is an interval) $\bigcup_{\tubes}T$ has volume $\approx 1$.


\subsection{The general case.}

In \cite{WangZahl2024}, Wang and the author considered the following variant of Definition \ref{defnWolffAxioms}.
\begin{definition}\label{defnFrostmanWolffAxioms}
Let $\tubes$ be a set of $\delta$ tubes in $\RR^n$. We say that $\tubes$ satisfies the \emph{Frostman Convex Wolff Axioms} if every convex set $W\subset\RR^n$ contains $O\big(|W|(\#\tubes)\big)$ tubes from $\tubes$. 

We say that $\tubes$ satisfies the \emph{Frostman Convex Wolff Axioms relative to the convex set} $V$ if the tubes in $\tubes$ are contained in $V$, and every convex set $W\subset V$ contains $O\big((|W|/|V|)(\#\tubes)\big)$ tubes from $\tubes$.
\end{definition}
Note that if $\#\tubes \sim \delta^{1-n}$, then $\tubes$ satisfies the Katz-Tao Convex Wolff Axioms (in the sense of Definition \ref{defnWolffAxioms}) if and only if it satisfies the Frostman Convex Wolff Axioms. If $\#\tubes<\!\!<\delta^{1-n}$ then it is impossible for $\tubes$ to satisfy the Frostman Convex Wolff Axioms, while if $\#\tubes>\!\!>\delta^{1-n}$ (and if the tubes are contained in the unit ball) then it is impossible for $\tubes$ to satisfy the Katz-Tao Convex Wolff Axioms.  Next, we introduce the following multi-scale version of Definition \ref{defnFrostmanWolffAxioms}, which generalizes the definition of stickiness (Definition \ref{defnSticky}).

\begin{definition}\label{frostmanEveryScale}
Let $\tubes$ be a set of $\delta$ tubes in $\RR^n$. We say that $\tubes$ satisfies the \emph{Frostman Convex Wolff Axioms at every scale} if for all $\delta\leq\rho\leq 1$, there exists a cover $\tubes_\rho$ of $\tubes$ with the property that for each $T_\rho\in\tubes_\rho$, the set $\tubes[T_\rho]$ satisfies the Frostman Convex Wolff Axioms relative to $T_\rho$.
\end{definition}

Definition \ref{frostmanEveryScale} generalizes Definition \ref{defnSticky}, since if $\#\tubes\sim\delta^{1-n}$ then $\tubes$ is sticky if and only if it satisfies the Frostman Convex Wolff Axioms at every scale. In \cite{WangZahl2024}, Wang and the author modified the arguments from \cite{WangZahl2022} to obtain the following generalization of Theorem \ref{stickyKakeyaThm}:

\begin{proposition}\label{FrostmanWolffEveryScale}
Let $\tubes$ be a set of $\delta$ tubes in $\RR^3$ that satisfy the Frostman Convex Wolff Axioms at every scale. For each $T\in\tubes$, let $Y(T)\subset T$ with $|Y(T)|\geq(\log 1/\delta)^{-1}|T|$. Then
\[
\Big|\bigcup_{\tubes}Y(T)\Big|\gtrapprox 1.
\]
\end{proposition}

In \cite{WangZahl2025}, Wang and the author combined Proposition \ref{FrostmanWolffEveryScale} with several types of induction on scale to prove Theorem \ref{WangZahlWolffAxiomThm}. As is often the case when employing induction, it is crucial to select the correct induction hypothesis. In hindsight, the condition that $\delta$ tubes point in $\delta$-separated directions is not preserved under the types of truncation and rescaling arguments that were used to prove Theorem \ref{WangZahlWolffAxiomThm}. Instead, the Katz-Tao and Frostman Convex Wolff Axioms proved to be more amenable to these types of transformations. In \cite{WangZahl2024}, Wang and the author introduced two assertions that quantify partial progress towards Theorem \ref{WangZahlWolffAxiomThm}. These definitions were later simplified by Guth \cite{Guth2025}. We present Guth's version below.

\begin{definition}\label{defnCDE}
Let $\sigma> 0$. 
\begin{itemize}
\item We say that \emph{Assertion $\cD(\sigma)$ is true} if the following holds:\\
Let $\tubes$ be a set of $\delta$ tubes that satisfy the Katz-Tao Convex Wolff Axioms. Then
\begin{equation}\label{defnCDEqn}
\Big|\bigcup_{\tubes}T\Big|\gtrapprox (\#\tubes)^{-\sigma}\Big(\sum_{T\in\tubes}T\Big) .
\end{equation}

\item  We say that \emph{Assertion $\cE(\sigma)$ is true} if the following holds:\\
Let $\tubes$ be a set of $\delta$ tubes that satisfy the Frostman Convex Wolff Axioms. Then
\begin{equation}\label{defnCEEstimate}
\Big|\bigcup_{\tubes}T\Big|\gtrapprox \delta^{2\sigma}\Big(\sum_{T\in\tubes}|T|\Big)^\sigma.
\end{equation}
\end{itemize}
\end{definition}

Assertion $\cD(1)$ is clearly true. To prove Theorem \ref{WangZahlWolffAxiomThm}, we must show that Assertion $\cD(\sigma)$ is true for every $\sigma>0$. We will prove this by iterating the following two results:

\begin{proposition}\label{CDImpliesCE}
For $\sigma>0$, Assertion $\cD(\sigma)$ implies Assertion $\cE(\sigma)$.
\end{proposition}

\begin{proposition}\label{CEImpliesImprovedCE}
There is a continuous function $g\colon (0,1]\to(0,1]$ with $g(\sigma)<\sigma$ so that for $\sigma>0$, Assertion $\cD(\sigma)$ plus Assertion $\cE(\sigma)$ implies Assertion $\cD(g(\sigma))$.
\end{proposition}


Proposition \ref{CEImpliesImprovedCE} is the more interesting (and challenging) of the two steps; we will discuss a few of the ideas used in the proof. Fix $\sigma>0$ and suppose that Assertions $\cD(\sigma)$ and $\cE(\sigma)$ are true.  Let $\tubes$ be a set of $\delta$ tubes that satisfy the Katz-Tao Convex Wolff Axioms. To simplify our discussion, we will consider the special case where $\#\tubes\sim\delta^{-2}$. Our goal is to prove that $\bigcup_{\tubes}T$ has volume substantially larger than what is guaranteed by the estimate \eqref{defnCDEqn}, i.e.~we wish to obtain an inequality of the form
\begin{equation}\label{desiredVolumeBdVignette}
\Big|\bigcup_{T\in\tubes}T\Big| \gtrapprox \delta^{2\sigma-\alpha},
\end{equation}
for some $\alpha>0$ that depends continuously on $\sigma$. Equivalently, we must show that 
\begin{equation}\label{desiredMuBd}
\mu_\tubes\lessapprox\delta^{-2\sigma+\alpha}.
\end{equation}

A second way to obtain  \eqref{desiredVolumeBdVignette} is to show there exists a scale $\tau >\!\!> \delta$ (here $\tau >\!\!> \delta$ means that $\tau$ is substantially larger than $\delta$) such that the union $\bigcup_{\tubes} T$ has larger than expected density at scale $\tau$; we will return to this later in the argument. 

Recall that if $\tubes$ is sticky, then for each scale $\delta<\rho<1$ it is possible to find a set $\tubes_\rho$ consisting of about $\rho^{-2}$ essentially distinct $\rho$ tubes (we say that two $\rho$ tubes $T,T'$ are \emph{essentially distinct} if $|T\cap T'|\leq (1/2)|T|$), each of which contain about $(\delta/\rho)^2$ tubes from $\tubes$. In this case we can apply the sticky Kakeya theorem (or more precisely, Proposition \ref{FrostmanWolffEveryScale}) to immediately obtain \eqref{desiredVolumeBdVignette}.

Let us instead consider an opposite extreme, where for every scale $\delta<\rho<1$, every cover $\tubes_\rho$ of $\tubes$ contains far more than $\rho^{-2}$ tubes (i.e.~every cover contains at least $\rho^{-2-\nu}$ tubes, and each set $\tubes[T_\rho]$ contains at most $\delta^\nu (\delta/\rho)^{-2}$ tubes, for some small quantity $\nu>0$).

\medskip

\noindent {\bf A fine-scale estimate.}\\ 
Fix a scale $\rho\in[\delta,1]$. For each $T_\rho\in\tubes_\rho$, the $\delta/\rho$ tubes in $\tubes^{T_\rho}$ satisfy the Katz-Tao Convex Wolff Axioms. Applying Assertion $\cD(\sigma,\omega)$, we conclude that 
\begin{equation}\label{muFineDefn}
\mu_{\operatorname{fine}}:=\sup_{T_\rho\in\tubes_\rho}\mu_{\tubes[T_\rho]}\lessapprox \sup_{T_\rho\in\tubes_\rho}(\#\tubes[T_\rho])^{\sigma}\leq \delta^{\nu\sigma} (\delta/\rho)^{-2\sigma}.
\end{equation}

Inequality \eqref{muFineDefn} bounds the typical intersection multiplicity of the $\delta$ tubes inside a common $\rho$ tube. Next, we define the quantity $\mu_{\operatorname{coarse}}$ as follows. For a typical point $x\in\bigcup_{\tubes}T$, there are about $\mu_{\operatorname{coarse}}$ distinct $\rho$ tubes $T_\rho\in\tubes_\rho$ with the property that $x\in\bigcup_{\tubes[T_\rho]}T$. With this definition, we have
\begin{equation}\label{mainMuBound}
\mu_{\tubes}\lesssim \mu_{\operatorname{fine}}\mu_{\operatorname{coarse}}.
\end{equation}

In the past, researchers considered a less efficient variant of \eqref{mainMuBound} of the form $ \mu \lesssim \mu_{\operatorname{fine}}  \mu_{\tubes_\rho}$, where $\mu_{\tubes_\rho}$ is the typical intersection multiplicity of the tubes in $\tubes_\rho$. Our use of the more efficient estimate \eqref{mainMuBound} is a key new ingredient in the proof.

\medskip

\begin{figure}[h!]
 \includegraphics[width=.35\linewidth]{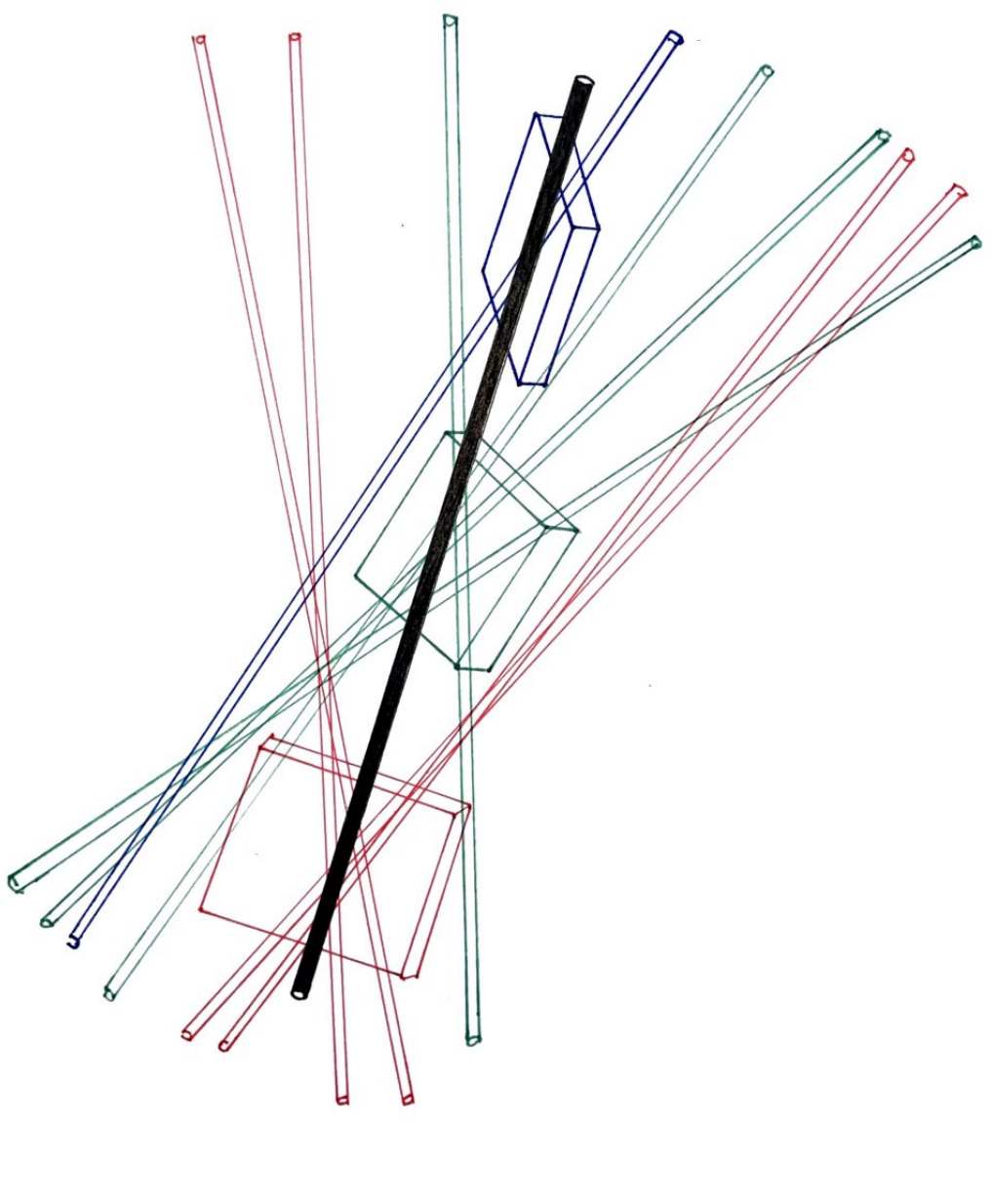}\hfill
 \includegraphics[width=.35\linewidth]{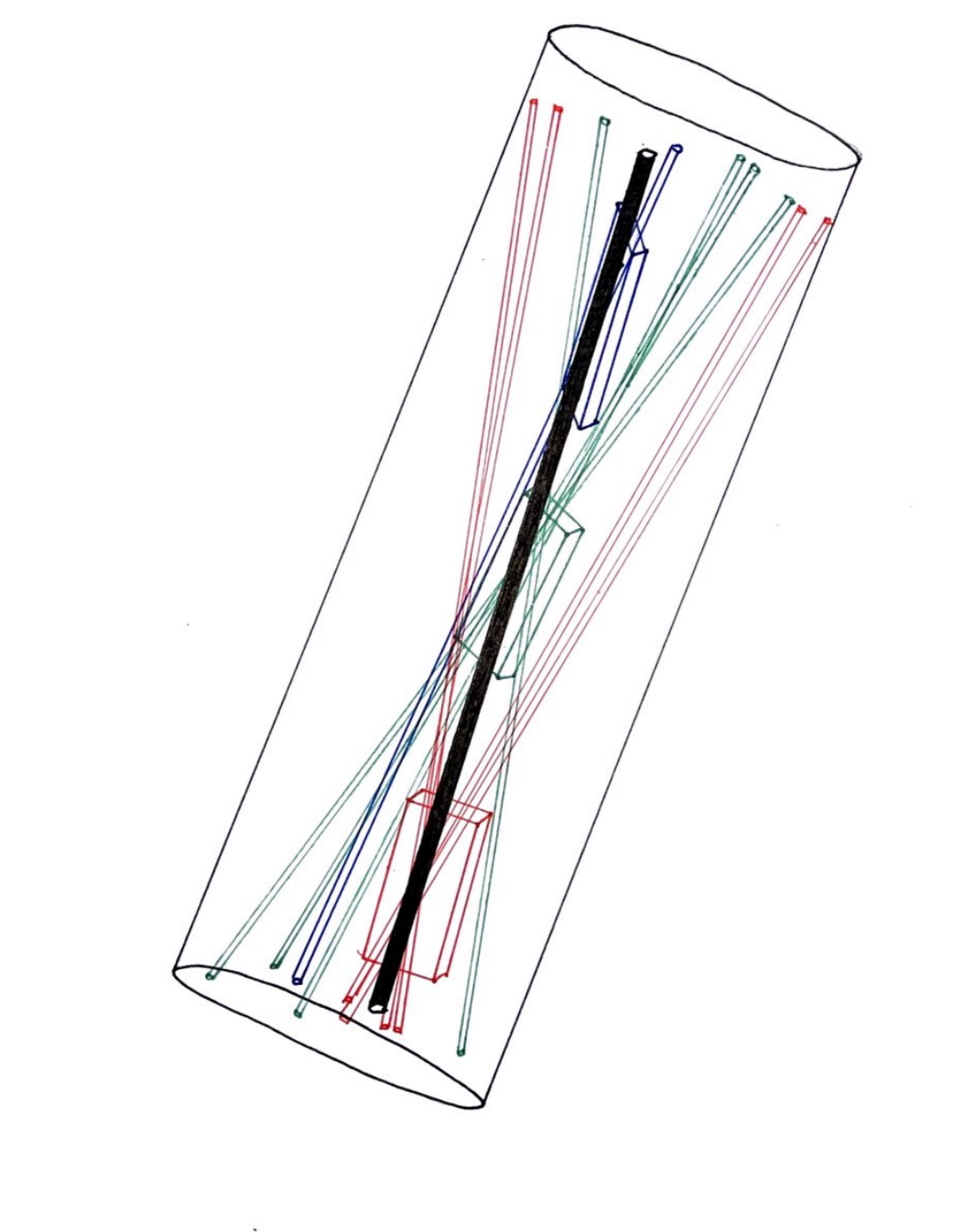}
\caption{Left: The set of tubes $\tubes^{T_\rho}$ and the corresponding grains. For clarity, we have only drawn the grains and tubes that intersect the black tube (and even most of these have been omitted; the set of red tubes passing through the red grain fill out a substantial portion of the red grain, and similarly for the other grains); the situation is similar for each tube in $\tubes^{T_\rho}$.\\
Right: The image of Figure \ref{grainsDecompHeartOfMatter} (left) after undoing the anisotropic rescaling that sent the tubes in $\tubes[T_\rho]$ to the tubes in $\tubes^{T_\rho}$. The dimensions of each grain have changed from $\delta/\rho\times c\times c$ to $\delta\times\rho c\times c$. }
\label{grainsDecompHeartOfMatter}
\end{figure}

\noindent {\bf A coarse-scale estimate.}\\
Our goal in this section is to establish the estimate 
\begin{equation}\label{desiredBoundMuCoarse}
\mu_{\operatorname{coarse}}\lessapprox \rho^{-2\sigma}.
\end{equation}
Combining \eqref{muFineDefn} and \eqref{desiredBoundMuCoarse}, we will obtain \eqref{desiredMuBd} with $\alpha = \nu\sigma$.  Fix a tube $T_\rho\in\tubes_\rho$. An argument due to Guth \cite{Guth2016b} says that the $\delta/\rho$ tubes in $\tubes^{T_\rho}$ arrange themselves into grains, i.e.~rectangular prisms; in this case these prisms will have dimensions $\delta/\rho\times c\times c$ for some $c>\!\!>\delta/\rho$. This means that we can cover $\bigcup_{\tubes^{T_\rho}}T^{T_\rho}$ by a set of (mostly) disjoint rectangular prisms of dimensions $\delta/\rho \times c\times c$; see Figure \ref{grainsDecompHeartOfMatter} (left).

Undoing the anisotropic rescaling that sent the tubes in $\tubes[T_\rho]$ to the tubes in $\tubes^{T_\rho}$, we have that $\bigcup_{\tubes[T_\rho]}T$ can be covered by a set of (mostly) disjoint rectangular prisms of dimensions $\delta \times \rho c \times c$; see Figure \ref{grainsDecompHeartOfMatter} (right). The same statement is true for each $T_\rho\in\tubes_\rho$. Let $\mathcal{G}$ denote the set of all such $\delta \times \rho c \times c$ prisms, from all $\rho$ tubes in $\tubes_\rho$. In order to bound $\mu_{\operatorname{coarse}}$, it suffices to bound the typical intersection multiplicity of the prisms in $\mathcal{G}$.

\medskip

Each $\delta \times \rho c \times c$ prism in $\mathcal{G}$ has an associated tangent plane, which is well-defined up to accuracy $\delta/(\rho c)$. If a typical pair of prisms intersect transversely, in the sense that their normal vectors make angle much larger than $\delta/(\rho c)$, then a straightforward $L^2$ argument shows that $\bigcup_{\tubes} T$ has larger than expected density at some scale $\tau >\!\!> \delta$ (the scale $\tau$ depends on the typical angle of intersection between prisms in $\mathcal{G}$), and as discussed above, we immediately obtain \eqref{desiredVolumeBdVignette}.

\begin{figure}[h!]
 \includegraphics[width=.3\linewidth]{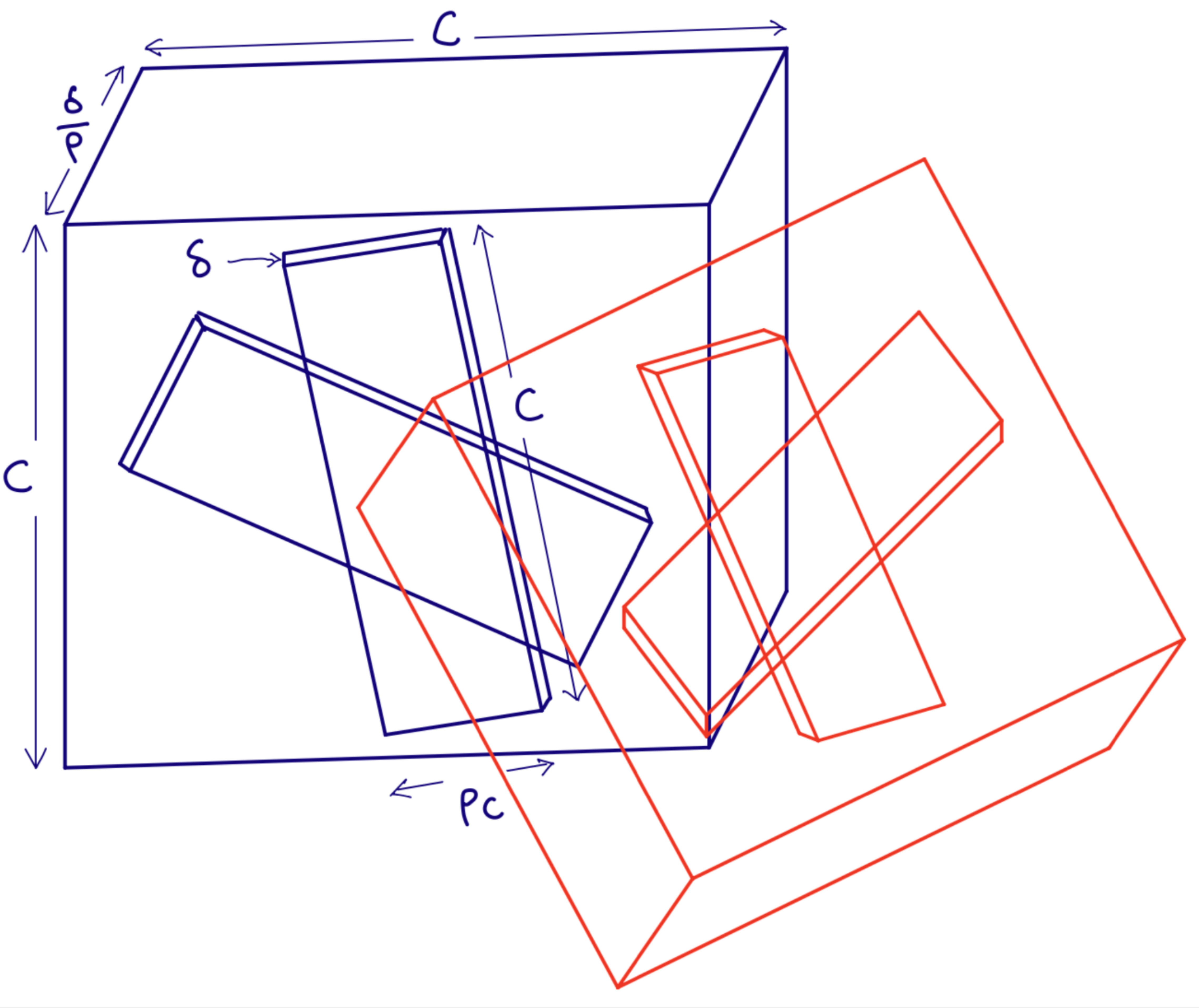}\hfill
 \includegraphics[width=.3\linewidth]{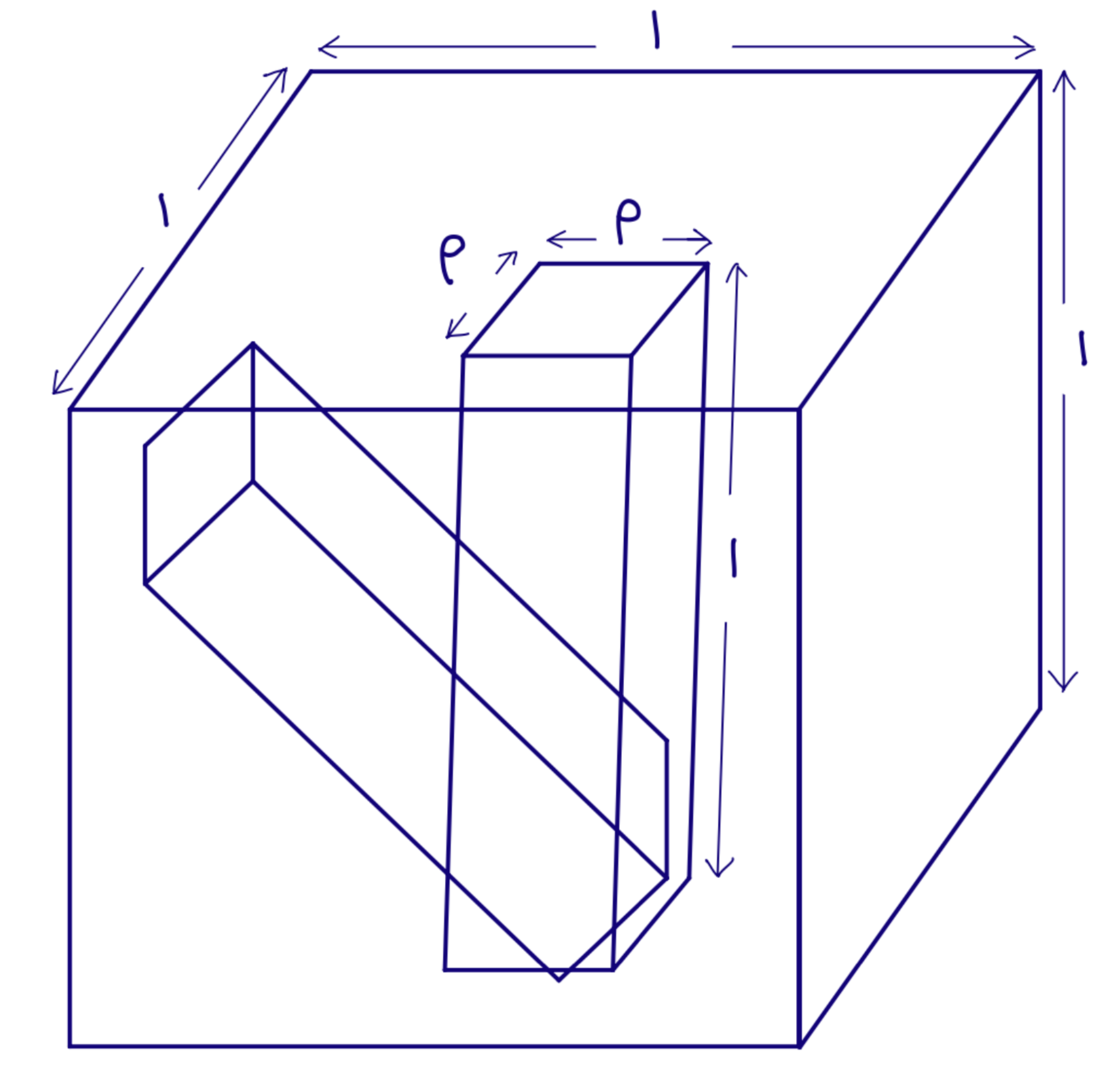}
\caption{Left: two sets of $\delta\times \rho c\times c$ prisms from the partition of $\mathcal{G}$ (blue and red, respectively), and the associated $\delta/\rho \times c \times c$ prisms ${\color{blue}\square}$ and ${\color{red}\square}$ that contain them.\\
Right: The anisotropic rescaling that maps the blue $\delta/\rho \times c \times c$ prism ${\color{blue}\square}$ to the unit cube maps each blue $\delta\times \rho c \times c$ prism to a $\rho\times\rho\times 1$ prism (this is comparable to a $\rho$ tube).}
\label{rescaledPrismsHeartOfMatter}
\end{figure}

Suppose instead that whenever two prisms $G,G'\in\mathcal{G}$ intersect, their corresponding tangent planes agree up to accuracy $\delta/(\rho c)$. We will call this the tangential case. This means that for each point $x$, the set of prisms from $\mathcal{G}$ containing $x$ are contained in a common prism of dimensions roughly $\delta/\rho\times c \times c$. Thus we can partition $\mathcal{G}$ into sets, $\mathcal{G} = \bigsqcup \mathcal{G}_i$, with the property that if two prisms intersect then they are contained in a common set, and the $\delta \times \rho c \times c$ prisms in each set $\mathcal{G}_i$ are contained in a common prism $\square_i$ of dimensions roughly $\delta/\rho \times c \times c$; see Figure \ref{rescaledPrismsHeartOfMatter} (left). This means that in order to bound $\mu_{\operatorname{coarse}}$, it suffices to bound the intersection multiplicity of the prisms within each set $\mathcal{G}_i$.

Fix a set $\mathcal{G}_i$, and let $\square_i$ be the associated $\delta/\rho\times c\times c$ prism. The image of each $G\in\mathcal{G}_i$ under the anisotropic rescaling sending $\square_i$ to the unit cube will be a prism of dimensions roughly $\rho\times\rho\times 1$ (see Figure \ref{rescaledPrismsHeartOfMatter} (right)); this is (morally) the same as a $\rho$ tube, i.e.~we can pretend that this set of prisms is actually a set of $\rho$ tubes; we will call this set $\tilde\tubes$. Our task of estimating $\mu_{\operatorname{coarse}}$ now reduces to estimating the typical intersection multiplicity of the tubes in  $\tilde\tubes$; our goal is to establish the bound
\begin{equation}\label{desiredBoundMuCoarseTildeTubes}
\mu_{\tilde\tubes}\lessapprox \rho^{-2\sigma}.
\end{equation}

We now consider the special case where the tubes in $\tilde\tubes$ satisfy the Katz-Tao Convex Wolff Axioms. If this is the case, then by applying Assertion $\cD(\sigma)$ we conclude that
\begin{equation}\label{muTildeTubes}
\mu_{\tilde\tubes}\lessapprox (\#\tilde\tubes)^{\sigma}\lesssim \rho^{-2\sigma},
\end{equation}
and we are done; for the second inequality we used the fact that every set of $\rho$ tubes satisfying the Katz-Tao Convex Wolff axioms has size $O(\rho^{-2})$. An important observation in the above argument is that after a suitable localization and rescaling, problem of understanding the intersection multiplicity of arrangements of grains can be converted into a problem of understanding the intersection multiplicity of arrangements of tubes; the latter can be studied using Assertion $\cD$ (this is a form of induction).

Finally, we will briefly discuss what happens when the tubes in $\tilde\tubes$ do not satisfy the Katz-Tao Convex Wolff Axioms. A key new idea is a structure theorem that finds a cover $\mathcal{W}$ of $\tilde\tubes$ with the property that $\mathcal{W}$ satisfies (a suitable analogue of) the Katz-Tao Convex Wolff Axioms, and for each $W\in \mathcal{W}$, the set $\tilde\tubes[W]=\{\tilde T\in\tilde\tubes\colon \tilde T\subset W\}$ 
satisfies the Frostman Convex Wolff Axioms. If the sets in $\mathcal{W}$ are ``thin,'' in the sense that their thinnest dimension is comparable to $\delta$, then we use an argument similar to the one described above with the sets in $\mathcal{W}$ in place of our previous grains. Otherwise, we apply (a suitably rescaled analogue of) Assertion $\cE(\sigma)$ to each set of tubes $\tilde\tubes[W],$ $W\in\mathcal{W}$ to conclude that $\bigcup_{\tubes} T$ has larger than expected density at a scale $\tau >\!\!> \delta$ (the choice of $\tau$ will depend on the thickness of the sets in $\mathcal{W}$), and as discussed above, we immediately obtain \eqref{desiredVolumeBdVignette}.

In the above proof, a few key new ideas stand out. First, we replace the hypothesis that tubes point in $\delta$-separated directions with the non-concentration conditions ``Katz-Tao Convex Wolff Axioms'' and ``Frostman Convex Wolff Axioms.'' These non-concentration conditions are carefully chosen to match the output of our structure theorem; after applying the structure theorem to an arbitrary collection of convex sets, we obtain a two-scale decomposition that satisfies the Katz-Tao Convex Wolff axioms at the coarse scale, and the Frostman Convex Wolff Axioms at the fine scale. Second, we use the efficient Inequality \eqref{mainMuBound} to combine estimates across different scales. Third, we apply induction on scale (specifically Assertion $\cD$) to control the interactions between (anisotropic) grains coming from different (coarse scale) tubes.

\subsection{Future directions.}
There appear to be several obstructions preventing the arguments in \cite{WangZahl2022, WangZahl2024, WangZahl2025} from being adapted to solve the Kakeya maximal function conjecture in $\RR^3$. Chief among them is the following issue: at many points in the argument, we must consider a set $\tubes$ of $\delta$ tubes, along with a shading $Y(T)\subset T$ satisfying $|Y(T)|\geq\lambda|T|$. In a small neighbourhood of a typical point $x\in \bigcup_\tubes Y(T)$, we apply a variant of Cordoba's two-dimensional Kakeya estimate to conclude that $\bigcup_\tubes Y(T)$ fills out a sizeable fraction of a $\delta\times a \times b$ rectangular prism, for some $\delta<\!\!<a\leq b\leq 1$. Unfortunately, Cordoba-type arguments typically show that $\bigcup_\tubes Y(T)$ fills out a $\lambda^2$-density fraction of this prism, rather than a $\lambda$-density fraction. After various rescaling arguments, this prism becomes a $\rho$ tube (for some $\delta\leq\rho\leq 1$); we then apply induction on scale to this set of $\rho$ tubes, with their associated $\geq\lambda^2$ density shading. If this induction step occurs many times, then the original $\lambda$-density shading becomes a $\lambda^{2^N}$ density shading, for some very large $N$. An interesting open problem is whether the arguments in \cite{WangZahl2022, WangZahl2024, WangZahl2025} can be performed in a more efficient manner (perhaps with stronger hypotheses about the shadings $Y(T)$) to prevent the losses described above. A second interesting question is to what extent the ideas from \cite{WangZahl2025} can be used to understand the Fourier restriction problem in $\RR^3$.


\section{Kakeya in higher dimensions.}

\subsection{Varieties containing many lines, and the polynomial Wolff axioms.}\label{manyLinesPolyWolffSection}
A plane in $\RR^3$ contains a two-dimensional family of lines. Similarly, the slab $[0,1]^2\times[-\delta,\delta]$ contains a set of $\delta$ tubes $\tubes$ of cardinality roughly $\delta^{-2}$, with $\big|\bigcup_\tubes T\big|\sim\delta$. This is not an interesting near-miss to the Kakeya set conjecture, because $\tubes$ (badly!) fails to satisfy the Convex Wolff Axioms. It is straightforward to show that a plane is the only (irreducible) algebraic variety in $\RR^3$ containing a two-dimensional family of lines, and thus the thin neighbourhood of an algebraic variety in $\RR^3$ is never an interesting near-miss to the Kakeya set conjecture.

In higher dimensions, however, the situation is more problematic. Define
\begin{equation}\label{badQuadraticHypersurface}
Z = \{(a,b,c,d)\in\RR^4\colon ad-bc=1\},\quad E = N_{\delta}(Z)\cap \overline{B(0,1)}.
\end{equation}
$Z$ is a quadric hypersurface that contains a 3-dimensional family of lines, and $E$ contains roughly $\delta^{-3}$ essentially distinct $\delta$ tubes. For each $x\in E$, there are approximately $\delta^{-1}$ distinct $\delta$ tubes that contain $x$; their union forms the $\delta$-neighbourhood of a cone with vertex $x$. The set $E$ contains about $\delta^{-3}$ essentially distinct tubes, but these tubes do not satisfy the Convex Wolff Axioms. To see this, fix a point $x\in E$, and let $E_x$ be the union of $\delta$-tubes that contain $x$ and are contained in $E$; this forms a thickened cone. Let $R$ be a $\delta\times\delta\times\delta^{1/2}\times 1$ rectangular prism that contains $x$ and is tangent to this cone; after a rotation and translation, we may suppose that $R = [0,\delta]\times[0,\delta]\times[0,\delta^{1/2}]\times [0,1]$, and $x = (0,0,0,0)$. $R$ contains about $\delta^{-1/2}$ essentially distinct tubes that contain $x$ and are contained in $E$. After replacing $E$ by its $\delta$ neighbourhood, we can verify that each of the $\delta$-separated points $x_j = (0, 0, \delta j, 0),$ $j=1,\ldots,\delta^{-1/2}$ are also contained in $E$, and for each of these points $x_j$, there are about $\delta^{-1/2}$ essentially distinct $\delta$ tubes that contain $x_j$ and are contained in $E$. This gives us $\delta^{-1/2}$ families, each with $\delta^{-1/2}$ tubes, for a total of $\delta^{-1}$ tubes (after a harmless refinement, these tubes are essentially distinct) contained inside a prism of volume $\delta^{5/2}$. In summary, the above computation shows that the density of tubes inside $R$ is too large by a multiplicative factor of $\delta^{-1/2}$.

Let $\tubes'$ be a set of roughly $\delta^{-5/2}$ tubes obtained by randomly selecting each tube contained in $E$ with probability $\delta^{1/2}$. For a typical $x\in E'$, the set of tubes from $\tubes'$ containing $x$ are still contained in the aforementioned thickened cone, but now (with a small number of exceptions), they point in $\delta^{1/2}$-separated directions. We can verify that the tubes in $\tubes'$ satisfy the Convex Wolff Axioms. 

Finally, it is straightforward to construct a set $\tubes$ consisting of roughly $\delta^{-1/2}$ translated or rotated copies of $\tubes'$, so that $\tubes$ has cardinality roughly $\delta^{-3}$, satisfies the Convex Wolff Axioms, and has volume $\big|\bigcup_\tubes T\big|\sim\delta^{1/2}$. Such a set $\tubes$ is a near-miss to the Kakeya set conjecture, and in particular Theorem \ref{WangZahlWolffAxiomThm} is false in dimension $n\geq 4$. In higher dimensions the situation becomes even worse, as the number and complexity of infinitely ruled surfaces increases. 

While the above example shows that the analogue of Theorem \ref{WangZahlWolffAxiomThm} is false in dimension $n\geq 4$, it is possible that the statement can be salvaged if we broaden the non-concentration condition to also exclude clustering inside sets such as \eqref{badQuadraticHypersurface}. 

\begin{definition}\label{polyWolffAxiomsDefn}
Let $\tubes$ be a set of $\delta$ tubes in $\RR^n$. We say that $\tubes$ satisfies the \emph{Polynomial Wolff Axioms} if for every semi-algebraic set $S$ and every $\delta\leq\lambda\leq 1$, we have
\[
\#\{T\in\tubes\colon |T\cap S|\geq\lambda\}\lesssim \lambda^{-n}|S|\delta^{1-n},
\]
where the implicit constant may depend on $n$ and the complexity of $S$ (i.e.~the number and degrees of the polynomials used to define $S$).
\end{definition}

\begin{conj}\label{KakeyaForPolynomialWolff}
(A): Let $\tubes$ be a set of $\delta$ tubes in $\RR^n$ that satisfy the Polynomial Wolff Axioms. Then
\[
\Big|\bigcup_{\tubes}T\Big| \gtrapprox \sum_\tubes |T|.
\]

\medskip
(B): For each $T\in\tubes$, let $Y(T)\subset T$ with $|Y(T)|\geq (\log 1/\delta)^{-1} |T|$. Then
\[
\Big|\bigcup_{\tubes}T\Big| \gtrapprox \sum_\tubes |T|.
\]

\medskip
(C):
\[
\Big\Vert \sum_{T\in\tubes}\chi_T\Big\Vert_{\frac{n}{n-1}}\lessapprox 1.
\] 
\end{conj}
Every set of tubes pointing in $\delta$-separated directions satisfies the polynomial Wolff axioms. This was proved for $n=3$ by Guth \cite{Guth2016}; $n=4$ by the author \cite{Zahl2018}, and for all $n$ by Katz and Rogers \cite{KatzRogers2018}. Thus versions (A) and (B) of Conjecture \ref{KakeyaForPolynomialWolff} imply the Minkowski and Hausdorff dimension version of the Kakeya set conjecture, while version (C) implies the Kakeya maximal function conjecture.
The remainder of this section will survey progress on the Kakeya conjecture for families of tubes satisfying the Polynomial Wolff Axioms.


\subsection{Multilinear Kakeya and the Bourgain-Guth multilinear $\to$ linear argument.}
Recall that in \cite{KatzLabaTao2000}, Katz, \L{}aba, and Tao proved that if $K$ is a (hypothetical) Besicovitch set in $\RR^3$ with upper Minkowski dimension close to $5/2$, then the discretization of $K$ at scale $\delta$ leads to a set $\tubes$ of $\delta$ tubes that are plany. The Bennett-Carbery-Tao multilinear Kakeya theorem \cite{BennettCarberyTao2006} asserts that this is true for every set of tubes for which $\big|\bigcup_\tubes T\big|$ is substantially smaller than $\sum_\tubes|T|$.

There are now three substantially different proofs of the multilinear Kakeya theorem. The original proof by Bennett, Carbery, and Tao used methods related to heat flow; an endpoint estimate by Guth \cite{Guth2010} used tools from algebraic topology, which later inspired the proof of the joints theorem \cite{GuthKatz2010} and subsequent developments in incidence geometry \cite{GuthKatz2015}; the third proof used the Loomis–Whitney inequality and induction on scale \cite{Guth2015}. The version stated below is due to Carberry and Valdimarsson \cite{CarberyValdimarsson2012}, whose proof is based on the second technique described above.

\begin{theorem}[Multilinear Kakeya theorem]\label{multilinearKakeyaThm}
Let $\tubes$ be a set of $\delta$ tubes in $\RR^n$ and let $2\leq k\leq n$. Then
\begin{equation}
\Big\Vert\Big(\sum_{T_1,\ldots, T_k \in \tubes}\chi_{T_1}\cdots\chi_{T_n}|v_1\wedge\ldots\wedge v_k|\Big)^{\frac{1}{k}}\Big\Vert_{\frac{k}{k-1}}
\leq C(n)\Big(\frac{1}{\delta}\Big)^{\frac{n}{k}-1}\sum_{T\in\tubes}|T|.
\end{equation}
In the above expression $v_i=\operatorname{dir}(T_i)$ is the direction of the tube $T_i$, and $|v_1\wedge\ldots\wedge v_k|$ is the $k$-dimensional volume of the parallelepiped spanned by $v_1,\ldots,v_k$.
\end{theorem}

Specializing to the case where $\#\tubes=\delta^{1-n}$, we conclude that if $\Big|\bigcup_{\tubes}T\Big| \leq \delta^{n-k}\theta$ for some $0<\theta<1$, then for a typical point $x\in\bigcup_\tubes T$ and a typical tuple $(T_1,\ldots,T_k)$ of tubes containing $x$, we must have that $T_1,\ldots,T_k$ are contained in the $\theta^{1/(k-1)}$-neighbourhood of an affine $(k-1)$-plane that contains $x$. After pigeonholing, we conclude that for a typical point $x\in\bigcup_\tubes T$ there exists an affine $(k-1)$-plane $H_x$ containing $x$ with the property that most of the tubes $\{T\in\tubes\colon x\in T\}$ are contained in the $\theta^{1/(k-1)}$-neighbourhood of $H_x$. We remark that this is consistent with the Heisenberg, $SL_2$, and quadric hypersurface examples (for $k=3,$ $k=3$, and $k=4$, respectively) discussed in Sections \ref{HeisenbergSL2ExampleSection} and \ref{manyLinesPolyWolffSection}.


In \cite{BourgainGuth2011}, Bourgain and Guth introduced a new method, which we will call the Bourgain-Guth multilinear $\to$ linear argument, to prove new bounds for the Fourier restriction conjecture. In the context of the Kakeya problem, the Bourgain-Guth multilinear $\to$ linear argument is as follows. Let $\tubes$ be a set of $\delta$ tubes in $\RR^n$ with $\#\tubes=\delta^{1-n}$, and let $2\leq k\leq n$ be an integer. Let $x\in\bigcup_\tubes T$ be a typical point, and let $(T_1,\ldots, T_k)$ be a typical $k$-tuple of tubes that contain $x$.

At one extreme, the tubes $T_1,\ldots,T_k$ point in (quantitatively) linearly independent directions, in the sense that their direction vectors satisfy $|v_1\wedge\cdots\wedge v_k|\sim 1$. If this is the case, then we can apply Theorem \ref{multilinearKakeyaThm} to conclude that $|\bigcup_\tubes T|\gtrsim \delta^{n-k}$. 

At the opposite extreme, suppose that the tubes $T_1,\ldots,T_k$ are contained in the $\delta$-neighbourhood of an affine $(k-1)$-plane. Since at most $O(\delta^{2-k})$ essentially distinct $\delta$ tubes can contain a common point and be contained in the $\delta$ neighbourhood of a $(k-1)$-plane, we conclude that  $\Big|\bigcup_\tubes T\Big|\gtrsim\delta^{k-2}$.

We therefore have the estimate
\begin{equation}\label{broadAndNarrowEstimate}
\Big|\bigcup_{\tubes}T\Big|\gtrsim \min\Big(\delta^{n-k},\ \delta^{k-2}\Big).
\end{equation}
The first estimate on the RHS of \eqref{broadAndNarrowEstimate} is called a \emph{multilinear estimate} or \emph{broad estimate}, and the second is called a \emph{narrow estimate}. By selecting $k=(n+2)/2$ for $n$ even and $k=(n+1)/2$ for $n$ odd, we conclude that 
$
\big|\bigcup_\tubes T\big|\gtrsim \delta^{\lfloor n/2\rfloor+1}.
$
To make the above argument rigorous, one must also consider intermediate scenarios between the two extremes described above. This can be done using induction on scale.

The argument described above uses very little information about Besicovitch sets---it exploits neither the direction-separated hypothesis, the Convex Wolff Axioms, nor the Polynomial Wolff Axioms. Indeed, it only makes use of the Frostman non-concentration condition that for $\delta\leq\rho\leq 1$, each $\rho$ tube contains $\lesssim (\rho/\delta)^{n-1}$ tubes from $\tubes$. In this sense, the above argument is sharp---see \cite{Zahl2023} for a precise formulation.

In the next two sections we will discuss improvements to both the broad and narrow estimates, which make use of additional information about Besicovitch sets.


\subsection{Broad estimates.}
Theorem \ref{multilinearKakeyaThm} makes no assumptions on the tubes in $\tubes$. At this level of generality the result is sharp. When applying Theorem \ref{multilinearKakeyaThm} in the multilinear $\to$ linear argument sketched above, however, we typically have additional information about $\tubes$. For example, the tubes in $\tubes$ might point in $\delta$-separated directions, and/or satisfy the Polynomial Wolff Axioms. With this additional constraint, it is possible to improve Theorem \ref{multilinearKakeyaThm} as follows.

\begin{theorem}\label{multiLinKakeyaDifferentDirections}
Let $\tubes$ be a set of $\delta$ tubes in $\RR^n$ that satisfy the polynomial Wolff axioms, and let $2\leq k\leq n$. Then
\begin{equation}\label{broadEstimate}
\Big\Vert\Big(\sum_{T_1,\ldots, T_k \in \tubes}\chi_{T_1}\cdots\chi_{T_n}|v_1\wedge\ldots\wedge v_k|^{\frac{k}{d}}\Big)^{\frac{1}{k}}\Big\Vert_{\frac{d}{d-1}}
\lessapprox \Big(\frac{1}{\delta}\Big)^{\frac{n}{d}-1}\Big(\sum_{T\in\tubes}|T|\Big)^{\frac{n(d-1)}{(n-1)d}},
\end{equation}
where
\[
d = \frac{n^2+k^2+n-k}{2n}.
\]
\end{theorem}
\begin{corollary}
Let $2\leq k\leq n$ and let  $\tubes$ be a set of $\delta$ tubes in $\RR^n$ pointing in $\delta$-separated directions, with $\#\tubes\sim\delta^{1-n}$.  Suppose that for a typical point $x\in\bigcup_{\tubes}T$ and a typical $k$-tuple $T_1,\ldots,T_k$ of tubes containing $x$, we have $|v_1\wedge\ldots\wedge v_k|\sim 1$. Then
\[
\Big|\bigcup_\tubes T\Big| \gtrapprox \delta^{n-d},\qquad d = \frac{n^2+k^2+n-k}{2n}.
\]

\end{corollary}

Theorem \ref{multiLinKakeyaDifferentDirections} was proved for $n=4,k=3$ by Guth and the author \cite{GuthZahl2018}. In higher dimensions Theorem \ref{multiLinKakeyaDifferentDirections} was proved independently and concurrently by the author \cite{Zahl2021} and by Hickman, Rogers, and Zhang \cite{HickmanRogersZhang2022} (the latter paper proves a $k$-broad rather than $k$-multilinear estimate). The results in \cite{Zahl2021,HickmanRogersZhang2022} are stated for sets of tubes pointing in $\delta$-separated directions, but they also apply to sets of tubes satisfying the Polynomial Wolff Axioms.

When the broad estimate from Theorem \ref{multiLinKakeyaDifferentDirections} is used as input to the Bourgain-Guth multilinear $\to$ linear argument, it gives an improved Kakeya maximal function bound for $n\geq 7$; it also gives new Hausdorff dimension bounds in certain dimensions $n\geq 7$ (for certain other dimensions $n\geq 7$ and for Minkowski dimension bounds, the best-known estimate comes from arguments based on Theorem \ref{sumDifferences}; see \cite{KatzTao2002b} for details).

The proof of Theorem \ref{multiLinKakeyaDifferentDirections} uses a technique developed by Guth \cite{Guth2016} in the context of the restriction conjecture. The idea is as follows. First, we consider the case $k = n-1$. Let $\tubes$ be a set of $\delta$ tubes. Using the Guth-Katz polynomial partitioning technique from \cite{GuthKatz2015}, we cover $\bigcup_\tubes T$ by a union of \emph{semi-algebraic grains}---a semi-algebraic grain is a set of the form $S = B(x,r)\cap N_\delta(Z)$, where $Z = Z(P)$ is the zero-set of a low degree polynomial. This is called a \emph{grains decomposition} of $\tubes$. The Guth-Katz polynomial partitioning technique controls the number of grains in this decomposition, as well as the number of tubes from $\tubes$ that intersect each grain. We then apply Theorem \ref{multilinearKakeyaThm} to the tubes that intersect each grain, and sum the resulting contributions. The estimate becomes stronger as the number of grains increases. We use the assumption that the tubes satisfy the Polynomial Wolff Axioms to ensure that the number of grains is suitably large. For $k < n-1$, we further decompose each grain in the above grains decomposition to obtain a multi-level grains decomposition of $\tubes$. The main new difficulty is to control the number of tubes that intersect each nested sequence of grains in this multi-level decomposition.


\subsection{Narrow estimates.}
In \cite{KatzZahl2021}, Katz and the author obtained the following narrow estimate in $\RR^4$:
\begin{theorem}\label{katzZahlPlanebrushThm}
Let $\tubes$ be a set of $\delta$ tubes in $\RR^4$ that point in $\delta$-separated directions, with $\#\tubes\sim \delta^{-3}$. Suppose that for each point $x\in\bigcup_{\tubes}T$, there is an affine 2-plane $H_x$ so that all tubes $\{T\in\tubes\colon x\in T\}$ are contained in the $O(\delta)$-neighbourhood of $H_x$. Then
\begin{equation}\label{lowerBdPlanyTubes}
\Big|\bigcup_{\tubes}T\Big|\gtrapprox \delta^{2/3}.
\end{equation}
\end{theorem}
When combined with the broad estimate from Theorem \ref{multiLinKakeyaDifferentDirections}, this yields new bounds on the dimension of Besicovitch sets in $\RR^4$; see \cite{Kulkarni2024} for related work in higher dimensions. 

The idea behind Theorem \ref{katzZahlPlanebrushThm} is as follows. Suppose that roughly $\mu$ tubes from $\tubes$ intersect each point $x\in\bigcup_{\tubes}T$. Fix a point $x$, and let $T_1\ldots,T_\mu$ be the tubes from $\tubes$ that intersect $x$. All of these tubes are contained in the $O(\delta)$ neighbourhood of a common affine 2-plane $H_x$. Let $\tubes'$ be the set of tubes that intersect $H_x$. Each of $T_1,\ldots,T_\mu$ intersect about $\delta^{-1}\mu$ tubes from $\tubes$, and thus we should expect that $\#\tubes'\gtrsim \mu(\delta^{-1}\mu) = \delta^{-1}\mu^2$. If the tubes in $\tubes'$ are mostly disjoint, then we have
\[
\Big| \bigcup_\tubes T\Big| \geq\Big|\bigcup_{\tubes'}T\Big|\gtrsim |T|\delta^{-1}\mu^{-2}\sim\delta^2\mu^{-2}.
\]
On the other hand, we have $\big| \bigcup_\tubes T\big|\sim \mu^{-1}$. Comparing these estimates we obtain \eqref{lowerBdPlanyTubes}.

Now suppose instead that the tubes in $\tubes'$ are not mostly disjoint. At the opposite extreme, for each point $y \in \bigcup_{\tubes'}T$ there are at least two tubes $T,T'\in\tubes'$ that intersect $y$ and make angle $\sim 1$ at this point of intersection. But recall that all of the tubes containing $y$ are contained in the $O(\delta$) neighbourhood of a common plane $H_y$, and hence $H_y$ is the plane spanned by the directions of $T$ and $T'$. This 2-plane has one-dimensional intersection with $H_x$. This in turn means that all of the tubes from $\tubes$ intersecting $y$ are contained in $\tubes'$. We thus obtain a decomposition
\[
\bigcup_\tubes T = \bigcup_{\tubes'} T\ \sqcup\ \bigcup_{\tubes\backslash \tubes'} T.
\]
We obtain a lower bound for the first term using a variant of Wolff's hairbrush argument, and we handle the second term by induction. 


\subsection{Future directions.}
It is an interesting question what improvements are possible for the broad and narrow estimates described above, and whether such improvements will play a role in the resolution of the Kakeya set conjecture in higher dimensions. The Bourgain-Guth multilinear $\to$ linear argument was developed in the context of the Fourier restriction conjecture. The improved broad estimate from Theorem \ref{multiLinKakeyaDifferentDirections} was applied to the Fourier restriction conjecture by Hickman and the author \cite{HickmanZahl2024} (see also \cite{Guth2016,Guth2018,HickmanRogers2019} for previous work in this direction), though this has since been superseded by different methods \cite{GuoWangZhang24, WangWu2024}. It seems plausible that improved broad and narrow estimates may play a role in future work on the restriction conjecture.

\section{Thanks.}
The author would like to thank Larry Guth and Hong Wang for comments and suggestions on an earlier draft of this manuscript. The author's research was supported by Discovery and Alliance grants from the Natural Sciences and Engineering Research Council of Canada; by the Fundamental Research Funds for the Central Universities (Project No.~100-63253272); and by the Nankai Zhide Foundation.
\bibliographystyle{plain}
\bibliography{ICM_bibligraphy}

\end{document}